\documentclass{article}

\usepackage{amsmath}
\usepackage{amsfonts}

\newtheorem{theo}{Theorem}

\newtheorem{prop}{Proposition}
\newtheorem{lemm}{Lemma}
\newtheorem{defn}{Definition}

\def\E{\mathbb{E}}
\def\0{{\bf 0}}
\def\Z{\mathbb{Z}}
\def\R{\mathbb{R}}

\renewcommand{\E}{\mathbb E \,}

\newcommand{\tod}{\stackrel{{\cal D}}{\longrightarrow}}

\newcommand{\eqco}{\setcounter{equation}{0}}
\newcommand{\thco}{\setcounter{theo}{0}}
\newcommand{\prco}{\setcounter{prop}{0}}
\newcommand{\laco}{\setcounter{lemm}{0}}
\newcommand{\coco}{\setcounter{coro}{0}}
\newcommand{\cjco}{\setcounter{conj}{0}}

\newcommand{\deco}{\setcounter{defn}{0}}
\newcommand{\allco}{\eqco  \thco \prco \laco \coco \cjco \deco}
\newcommand{\qed}{\hfill{\rule[-.2mm]{3mm}{3mm}}}

\newcommand{\Po}{{\cal P}}
\newcommand{\Q}{{\cal Q}}

\renewcommand{\P}{{\cal P}}

\newcommand{\X}{{\cal X}}

\newcommand{\NN}{{\cal N}}

\def\bdm{\begin{displaymath}}
\newcommand{\edm}{\end{displaymath}}
\def\benu{\begin{enumerate}}
\def\eenu{\end{enumerate}}
\def\beqn{\begin{equation}}
\def\eeqn{\end{equation}}
\def\be{\begin{equation}}
\def\ee{\end{equation}}
\def\bea{\begin{eqnarray}}
\def\eea{\end{eqnarray}}
\newcommand{\bean}{\begin{eqnarray*}}
\newcommand{\eean}{\end{eqnarray*}}
\newcommand{\bear}{\begin{eqnarray}}
\newcommand{\eear}{\end{eqnarray}}

\renewcommand{\epsilon}{\varepsilon}

\newcommand{\out}{{\rm out}}
\newcommand{\inn}{{\rm in}}

\begin{document}

\title{\bf
    Gaussian Fields and Random Packing}

\bigskip

\author{Yu. Baryshnikov and J. E. Yukich$^{1}$ \\
\\
{\normalsize{\em Bell Laboratories and Lehigh University}} }

\date{\today}
\maketitle

\footnotetext[1] {Rm 2C-323,
Bell Laboratories, Lucent Technologies, 600-700 Mountain Ave, Murray Hill,
NJ 079074:
{\texttt ymb@research.bell-labs.com} }

\footnotetext[2]{ Department of Mathematics, Lehigh University,
Bethlehem PA 18015, USA: {\texttt joseph.yukich@lehigh.edu} }

\footnotetext{$~^1$ Research supported in part by NSA grant
MDA904-01-1-0029 }

\begin{abstract} Consider sequential packing of unit volume balls
in a large cube, in any dimension and with Poisson input.  We show
after suitable rescaling that the spatial distribution of packed
balls tends to that of a Gaussian field in the thermodynamic
limit.  The results cover  related applied models, including
ballistic deposition and spatial birth-growth models.

\end{abstract}

\section{ Introduction}

The following prototype random packing model is known as the basic
Random Sequential Adsorption Model (RSA) for hard spheres on a
continuum surface. Open balls  $B_{1,n},B_{2,n}...,$ of unit
radius arrive sequentially and uniformly at random in the
$d$-dimensional cube $Q_n$ having volume $n$ and centered at the
origin. Let the first ball $B_{1,n}$ be {\em packed}, and
recursively for $i=2,3, \ldots$, let
 the $i$-th ball
$B_{i,n}$ be  packed iff $B_{i,n}$ does not overlap any ball in
$B_{1,n},...,B_{i-1,n}$ which has already been packed. If not
packed, the $i$-th ball is discarded.
 Given a positive integer $k$, let
$  N(\{B_{1,n},...,B_{k,n} \})$ be the number of balls  packed out
of the first $k$ arrivals. $N_{n,d}(k):= N(\{B_{1,n},...,B_{k,n}
\})$ are called random packing numbers.

Lattice packing is defined analogously to continuum packing, save
for the obvious constraint that centers of incoming balls are
constrained to lie on a lattice.

There is a vast literature involving versions of the  RSA model on
continuum and lattice substrates. There is a plethora of
experimental results and a notable dearth of mathematically
rigorous results, particularly in more than one dimension.  For
surveys of the extensive literature, see Evans \cite{E}, Senger et
al. \cite{SVS},
 Bartelt and Privman \cite{BP}, Adamczyk et al. \cite{ASZB}, Talbot et al.
\cite{T} and \cite{Pr}.

In addition to their fundamental role in adsorption modelling,
sequential packing models
 arise  in the study of hard core interactions in
physical and materials science,  spatial growth models in
crystallography and biology (Evans \cite{E}, sect. III,
Garcia-Ruiz et al. \cite{GLMS}), and in the study of polymer
reactions \cite{GLMS, Pr}. In modelling communication protocols
(Coffman et al. \cite{CFJP}), RSA is
 called {\em on-line packing}.

Consider input of size $\Po(\tau n)$, that is  the deposition
intensity is a Poisson random variable with parameter $\tau$.
Deposition intensity in the continuum is the average number of
particles arriving per unit volume, whereas in the lattice it is
the average number of arriving particles per lattice point. In
\cite{PY2}, the authors show that for both continuum and lattice
packing, the Poisson packing numbers $N_{n,d}(\Po(\tau n))$
satisfy a thermodynamic limit as well as a central limit theorem.

\begin{theo}
\label{main0} (LLN and CLT for  packing numbers \cite{PY2}) For
all $\tau \in (0,\infty)$ and all $d \geq 1$ there are positive
constants $\alpha_{d, \tau}$ and  $ \eta_{d,\tau}$ such that
\begin{equation}
\frac{ N_{n,d}(\Po(\tau n))} {n}  \to \alpha_{d, \tau} \ \ \ \
c.m.c.c. \label{LLN1}
\end{equation}
while
\begin{equation}
\frac{ N_{n,d}(\Po(\tau n)) - EN_{n,d}( \Po(\tau n))  } {n^{1/2} }
 \tod \NN(0, \eta_{{d,\tau}}^2 )
\label{CLT1}
\end{equation}
and
\begin{equation}
n^{-1} {\rm Var\,} N_{n,d}(\Po[\tau n ])  \to \eta_{d,\tau}^2.
\label{CLT2}
\end{equation}

\end{theo}

Here c.m.c.c. denotes convergence of means and complete
convergence, $\alpha_{d, \tau}$ is called the coverage function,
and $\NN(0, \sigma^2)$ denotes a normal random variable with mean
zero and variance $\sigma^2$.  The limits (\ref{LLN1}-\ref{CLT2})
also hold if the Poisson packing numbers $N_{n,d}(\Po(\tau n))$
are replaced by the finite input packing numbers $N_{n,d}([\tau
n])$ \cite{PY2}.

 Our goal here is to investigate the large $n$ {\em
distribution} of packed balls on the d-dimensional substrate
$Q_n$. Even in dimension $d = 1$  little is known about the
distribution of the point process on $\R$ induced by the packed
balls.  Our main result is that the basic RSA packing  process on
$Q_n \subset \R^d$  induces a point measure which, when suitably
rescaled, converges to that of a Gaussian field on $\R^d$ in the
thermodynamic limit. The CLT given by (\ref{CLT1}) of Theorem
\ref{main0} is a by-product. Our main result also holds for
variants of the basic packing process.

The packing process is defined on the infinite substrate $\R^d$ in
a natural way.  In this context, we represent the centers of the
incoming balls, together with their arrival times, as a point set
in $\R^d \times \R^+ $. Points of $\R^d \times \R^+$ are
generically denoted by $w:=(x,t_x)$, where $x \in \R^d, \ t_x \in
\R^+$. Given a locally finite point set $\X \subset \R^d \times
\R^+$, we let $\alpha(\X)$ denote the subset of $\X$ which is
accepted in the packing process.  We let $\pi(\X) \subset \R^d$
denote the projection of $\alpha(\X)$ onto $\R^d$. $\pi(\X)$ is
the point process on the substrate $\R^d$ formed by the accepted
balls.

 Let $\Po_{\tau}$ denote a rate one
homogeneous Poisson point process on $\R^d \times [0, \tau]$. In
the lattice setting,  ${\cal Q}$ denotes a collection of rate one
homogeneous Poisson point processes on $\R$ indexed by the lattice
points $\Z^d$ and embedded in a natural way into  $\R^d \times
\R^+ $. $\Po_{\tau}$ and ${\cal Q}$ henceforth represent the input
for the continuum and lattice packing processes, respectively.
The processes $\pi(\Po_{\tau})$ and $\pi(\cal Q)$ are thinned
Poisson point processes; we will be interested in their spatial
distribution.

For any set $A \subset \R^d$, let $\Po_{\tau, A}$ be the
restriction of $\Po_{\tau}$ to $A \times [0, \tau]$ and let ${\cal
Q}_A$ denote the restriction of ${\cal Q}$ to $A$.

For any Borel set $B \subset \R^d$ define the random fields
$$
\nu_{\tau}(B):= \sum_{x \in \pi (\Po_{\tau})} \delta_x(B)
$$
and
$$
\nu_{\tau,A}(B):= \sum_{x \in \pi (\Po_{\tau, A})} \delta_x(B);
$$
that is, $\nu_{\tau}$ is the random field on subsets of $\R^d$
induced by the  packing process on $\R^d$ whereas $\nu_{\tau, A}$
is the random field on subsets of $A$ induced by the packing
process on $A$. Since $\sup_B \nu_{\tau}(B)  = \infty$ and $\sup_B
\nu_{\tau, A}(B) \leq \text{volume}(A)$, we call $\nu_{\tau}$ and
$\nu_{\tau, A}$ the {\em infinite volume} and {\em finite volume}
packing measures, respectively. Edge effects show that in general
these two random fields do not coincide on subsets of $A$.

We similarly define random fields induced by the lattice packing
process as follows.  For any Borel set $B \subset \R^d$ define
$$
\mu(B):= \sum_{x \in \pi ({\cal  Q})} \delta_x(B)
$$
and
$$
\mu_A(B):= \sum_{x \in \pi ({\cal  Q}_A)} \delta_x(B).
$$

 Consider the rescaled infinite volume continuum  packing measures
\begin{equation} \label{measure}
\nu_{\tau, \lambda}(B):= \frac{\nu_{\tau} (\lambda B)- E
\nu_{\tau} (\lambda B)} {\sqrt{\lambda^d}}
\end{equation}
and the rescaled finite volume continuum packing measures
\begin{equation}\label{measure1}
\nu_{\tau, A, \lambda}(B):= \frac{\nu_{\tau, \lambda A} (\lambda
B)- E \nu_{\tau, \lambda A} (\lambda B)} {\sqrt{\lambda^d}}.
\end{equation}
The rescaled centered lattice packing measures $\mu_{\lambda}$ and
$\mu_{A, \lambda}$ are defined analogously.  All packing measures
$\mu$ (or $\nu$) considered here are spatially homogeneous in the
sense that the vectors $\langle \nu(B_1),...,\nu(B_k) \rangle$ and
$\langle\nu(B_1 + x),...,\nu(B_k + x) \rangle$ have the same
distribution for any $x \in \R^d$ and any Borel subsets
$B_1,...,B_k$.

In what follows all random measures are defined on the Borel
subsets of $\R^d$.
 Recall that a sequence of random fields $\mu_n, \ n \geq 1,$
converges in distribution to $\mu$ if and only if all finite
dimensional distributions $\langle \mu_n(B_1),...,\mu_n(B_k)
\rangle$ converge to $\langle\mu(B_1),...,\mu(B_k) \rangle$, where
$B_1,...,B_k$ are Borel subsets of $\R^d$. The following are our
main results.

\begin{theo} (Infinite volume packing measures converge to a Gaussian field) \label{main1}

(i) (Poisson input continuum packing) For all $\tau < \infty$,
$\nu_{\tau, \lambda}$ converges in distribution as $\lambda \to
\infty$ to a generalized Gaussian random field with covariance
kernel $K_{\nu, \tau}$ concentrated on the diagonal, that is
$$
K_{\nu, \tau} (x,y) :=C_{\nu, \tau} \delta(x-y),
$$
where $C_{\nu, \tau}$ is a constant depending on $\tau$.

(ii)(infinite input lattice packing)
 $\mu_{ \lambda}$
 converges  in distribution as $\lambda \to
\infty$ to a generalized Gaussian random field with covariance
kernel $K_{\mu}$ concentrated on the diagonal, that is
$$
K_{\mu}(x,y) := C_{\mu} \delta(x-y),
$$
where $C_{\mu}$ is a constant.
\end{theo}

The next theorem clearly extends Theorem \ref{main0}. To obtain
Theorem \ref{main0} we simply let $A$ be the unit cube centered at
the origin of $\R^d$.

\begin{theo} (Finite volume packing measures converge to a Gaussian field) \label{main2}
Let $A \subset \R^d$  have piecewise smooth boundary.  Then:

(i) (Poisson input continuum packing) For all $\tau < \infty$,
$\nu_{\tau, A, \lambda}$  converges  in distribution as $\lambda
\to \infty$ to a generalized Gaussian random field with covariance
kernel $K_{\nu, \tau}$.

(ii)(infinite input lattice packing)
 $\mu_{A, \lambda}$ converges  in distribution as $\lambda \to
\infty$ to a generalized Gaussian random field with covariance
kernel $K_{\mu}$.
\end{theo}

{\em Remarks.} (i) We will see in the sequel that
$$ C_{\nu, \tau} :=  \int_V [r_2(0,x)-r_1(0)r_1(x)] dx + r_1(0)
,
$$
where $r_1$ and $r_2$ are the one and two point correlation
functions for the spatial point process of packed points on
$\R^d$.  There is a similar expression for $C_{\mu}$. Neither
$K_{\nu, \tau}$ nor $K_{\mu}$ depend on the set $A$.

(ii)  By definition of weak convergence,
 Theorem \ref{main1} tells us that for any Borel sets
$B_1,..., B_m$ in $\R^d$ the $m$-vector with entries given by the
$\nu_{\tau, \lambda}$ measure of $B_1,..., B_m$ tends to a
Gaussian limit with  covariance matrix $C_{\nu, \tau}
\text{volume}(B_i \cap B_j)$.  Theorem \ref{main2} makes an
analogous statement for the $m$-vector with entries given by the
$\nu_{\tau, A, \lambda}$ measure of $B_1,..., B_m$.

(iii) Proving convergence to a Gaussian field for infinite input
continuum packing remains an open problem. In dimension $d = 1$,
Dvoretzky and Robbins \cite{DR} showed that the number of packed
balls asymptotically converges to a normal random variable.
Theorem \ref{main2}(ii) adds to results of Penrose \cite{P2}, who
shows that the number of packed balls in the lattice setting
satisfies a central limit theorem.

\vskip.5cm

{\em Acknowledgements.}  The authors are grateful for
conversations with Y. Suhov, whose questions helped inspire the
investigations in this paper.

\section{Related Models}
\allco

There are a multitude of variants of the basic RSA packing model
\cite{E,Pr}.  The approach taken here shows that the spatial
distribution of packed balls (particles) for many of these
variants converges to a Gaussian field. We discuss these variants
below.  Formal details may be found in \cite{PY2}.

\subsection{Random shapes and types}
\allco

The basic RSA model assumes that the incoming particle is a ball
of unit radius.  However, this assumption may be relaxed to allow
random shapes and random sets in $\R^d$, a widely considered model
\cite{BP}. More generally, arriving particles may have a random
type or spin, not necessarily representing shape or size. In one
dimension such a model is considered by Itoh and Shepp \cite{IS}.
Instead of considering the point process defined by the point set
of accepted particles, we may consider the point process
consisting of only those accepted particles of a particular kind
or type. Theorems \ref{main1} and \ref{main2} hold for the
measures induced by such point processes.

\subsection{Time dependent models}
\allco The basic RSA packing model can be generalized to include
the case in which a packed ball remains in place for a random
period of time at the end of which it is removed, i.e., desorbs.
This is a dynamic model, which among other things, describes the
reversible deposition of particles on substrates \cite{SVS, T}.
Assuming that the spatial locations and arrival times of particles
are given by a space-time Poisson process of unit intensity on
$\R^d \times [0, \tau]$, the contribution of a particle to the
point process of adsorbed points is determined not only by whether
it is accepted or not, but also by whether, if accepted, it
desorbs by time $\tau$.  The point measures given by the process
of adsorbed points are defined analogously to (\ref{measure},
\ref{measure1}) and in the large $\lambda$ limit satisfy the
Gaussian structure results of Theorems \ref{main1} and
\ref{main2}.

We can also extend the basic packing model to a generalized
version of the classical birth-growth model on $\R^d$ in which
cells are formed at random locations $x_i \in \R^d$ at times $t_i,
\ i = 1,2,...$ according to a unit intensity homogeneous
space-time Poisson point process.  When a new cell is formed, its
center $x_i$ is called its ``seed".  Once a seed is born, the cell
around it initially takes the form of a ball of possibly random
radius. If the initial radii of cells are zero a.s., then the
model is referred to as the Johnson-Mehl model (see Stoyan and
Stoyan \cite{SS}).
 The ball immediately generates a cell by growing radially
in all directions with a constant speed. Whenever one cell touches
another, it stops growing in that direction.
New seeds, and the  cells around them, form only in the uncovered
space in $\R^d$. The point measures given by the process of seeds
are defined analogously to (\ref{measure}, \ref{measure1}) and in
the large $\lambda$ limit satisfy the Gaussian structure results
of Theorems \ref{main1} and \ref{main2}. This adds to results of
Chiu and Quine \cite{CQ}, who only study the restriction to a
large window of
the infinite stationary birth-growth model on $\R^d$. Furthermore,
as in Penrose and Yukich \cite{PY2} (Theorem 2.1(b)), they
consider only the number of seeds generated and not their spatial
distribution.

\subsection{Ballistic Deposition Models}
\allco

The standard random particle deposition model considers random
size i.i.d. $(d+1)$-dimensional balls (``particles'') which rain
down sequentially at random onto a $d$ dimensional substrate of
volume $n$, and centered around the origin.

When a particle arrives on the existing agglomeration of deposited
balls,  the particle may slip and roll over existing particles,
undergoing displacements, stopping when it reaches a position of
lower height (the surface relaxation (SR) model \cite{BS}).  If a
particle reaches the surface $\R^d$, it is irreversibly fixed on
it; otherwise, the particle is removed from the system and the
next sequenced particle is considered. The rolling process does
not displace already deposited particles; there is no updating of
existing particles.  For $d = 1$, the model dates back to Solomon
\cite{So}. Senger et al. \cite{SVS} describe the many experimental
results.

The accepted particles all lie on the substrate and are
represented by points in $\R^d$.  The position of an accepted
particle is a translate of the original location in $\R^d$ above
which it originally comes in.   The point process of accepted
particles defines a point measure on $\R^d$.  The rescaled
measures are defined analogously to (\ref{measure}) and
(\ref{measure1}).  The methods described below may be easily
modified to show that in the large $\lambda$ limit, these measures
converge to a Gaussian random field. This thus generalizes the
results of \cite{PY2} (Theorem 2.2(b)), which only considers the
distribution of the number of accepted particles in the substrate.

\section{Auxiliary Results}
\allco

Throughout we let $W := V \times \R^+$, where $V := \R^d$. Recall
that points of $W$ are generically denoted by $w := (x,t_x)$,
where $x \in V, \ t_x \in \R^+.$ The spatial location of the
incoming balls, together with their arrival times, are represented
by points in $W$. The interaction range of arriving balls is just
the common diameter of the balls, that is to say equals two.
Throughout, $|| \cdot ||$ denotes the Euclidean norm on $\R^d$.

\subsection{Exponential Decay}

Let $\X \subset W$ be a locally finite point process.
 Following \cite{PY2} we make $\X$  into
the vertex set of an oriented graph by including an edge from
$w_1:=(x_1,t_{x_1})$ to $w_2:=(x_2,t_{x_2})$ whenever $||x_1-
x_2|| \leq 2$ and  $t_{x_1} \leq t_{x_2}$. Given $w \in \X$, let
$A_{\out}(w,\X)$ be the set of points (forward cone) in $\X$ that
can be reached from $w$ by a directed path  in this graph (along
with $w$ itself). Let $A_{\inn}(w,\X)$ be the set of points
(backward cone) in $\X$ from which the point $w$ can be reached by
a directed path in this graph (along with $w$ itself). Finally,
consider the ``causal cone''
$$
A_{\out, \inn}(w, \X ):= A_{\out}(w,\X) \cup A_{\inn}(w,\X).
$$
The next result justifies describing the sets $A_{\out, \inn}(w,
\Po_{\tau})$ and $A_{\out, \inn}(w, \Q)$ as ``cones''.

\begin{lemm}\label{cluster}
 Fix $\tau$ and let $\X$ be either
$\Po_\tau$ or $\Q$.
 There exist positive constants $\gamma:=\gamma(\tau)$, $\beta:=\beta(\tau),$
  and $\rho:=\rho(\tau)$,
 such that the causal cone associated with $w := (x,t_x) \in \X$
  belongs to the set
$$
C_R(w):= \{ (y,t_y): ||x-y|| \leq \beta |t_x - t_y| + R\}
$$
with probability at least  $1 - \rho \exp(-\gamma R)$.
\end{lemm}

\def\bbi{{\bold i}}

{\em Proof.} For the case that $\X$ is the finite input set
$\Po_\tau$, this result is just Lemma 4.2 of \cite{PY2}. For the
case that $\X$ is $\Q$,  we will use an ``invasion percolation''
argument similar to that in \cite{PY2}. The following argument
actually also holds if $\X$ is the infinite continuum input $\Po$.
 First, we discretize the
space $V$ by tiling it by lattice cubes of sufficiently large size
so that no ball can intersect non-adjacent cubes. More precisely,
we choose a full rank lattice $L\subset V$ (a sublattice of $\Z^d$
in the lattice case) with the elementary cube $C$,
$V=\cup_{\bbi\in L} C_\bbi$, where $C_\bbi:= C+\bbi$ and the size
of the cube $C$ is chosen as to satisfy $C+B\subset 2C$ (a further
restriction in the lattice case is that the boundary of $C$ does
not contain any point with integer coordinates, to avoid
ambiguity).

For each cube $C_\bbi, \ \bbi \in L,$ consider the arrival times
of the points whose centers belong to this cube; they form a
family $P(\bbi), \ \bbi\in L,$ of independent Poisson point
processes.

Now, fix a point $w=(x,t_x)$ and assume, without loss of
generality, that $x$ belongs to the cube centered at the origin.
Let us restrict attention to the {\it forward} cone and, to
simplify notation, assume that $t_x=0$, i.e, $w = (x, 0)$.  Call a
{\it path} $\pi := \pi_0,\pi_1,\ldots,\pi_p$ (of length
$p:=|\pi|$) a collection of points in $L$, such that for any two
consecutive points $\pi_k$ and $\pi_{k + 1}$, the cubes
$C_{\pi_k}$ and $C_{\pi_{k+1}}$ share at least one vertex. For a
path $\pi$ of length $p$, call the increasing sequence of  arrival
times $0=t_0<\ldots < t_p$ {\it admissible} if $t_k\in P(\pi_k)$
for all $0 \leq k \leq p$. Note that any path $\pi$ admits many
admissible sequences. Given a path $\pi$ of length $p$, let
$T(\pi)$ be the infimum of the terminal values $t_p$.



The relations above imply that if the forward cone of $(x,0)$ does
not belong to $C_R(w)$, then for some $\bbi$ there exists a path
$\pi$ from $\0$ to $\bbi$ such that $|\bbi | > \beta T(\pi) + R$,
that is,
$$
T(\pi)<(|\bbi|-R)/\beta.
$$

The probability that such a path exists is majorized by
the sum over all possible paths $\pi$, starting at the origin, of
the probability
$$
P[T(\pi)<(|\pi_{|\pi|}|-R)/\beta] \leq P[T(\pi)<(k|\pi|-R)/\beta]
$$
(we use here the obvious fact that the distance between the last
point of a path of length $p$ and the origin is at most $kp$,
where $k$ is a constant depending only on the lattice $L$ and the
dimension $d$).

The last piece we need before starting calculations is the fact
that for a path $\pi$ of length $p$, $T(\pi)$ is simply the sum of
$p$ \ i.i.d. \ exponential variables $Y$.
Choose $\theta$ large enough so that $E \exp(-\theta Y) \leq
3^{-(d + 2)}$ and choose $\beta$ large enough so that $e^{\theta
k/\beta} \leq 3$. Let $a := (kp - R)/\beta$.  Then
$$
P(\pi):=P[T(\pi) < a] \leq e^{\theta a} (E[e^{-\theta Y}])^p.
$$
By definition of $a$ and the choice of $\theta$ and $\beta$ we
have that the above is bounded by
$$
e^{\theta kp/\beta} 3^{-p(d + 2)} \leq 3^{-p(d + 1)}.
$$
Since the number of paths starting at $\0$ of length $p$ is
bounded by $3^{dp}$, the probability that the forward cone of
$(x,0)$ does not belong to $C_R(w)$ is bounded above by the sum of
the terms $P(\pi)$ over all possible paths $\pi$ starting at $\0$
and of length at least $R/k$.  This probability is thus bounded by
$$
\sum_{p \geq R/k} 3^{dp} P[T(\pi) < a] \leq \sum_{p \geq R/k}
3^{dp} 3^{-p(d + 1)},
$$
which gives the desired exponential decay in $R$.

Since  the case of the backward cone may be handled by using a
``time-reversal'', this completes the proof of Lemma 3.1. \qed

 Let $\X \subset W$ be a point set.  In keeping with the
terminology of statistical mechanics, we write $\sigma(w,\X)=1$,
if the ball centered at $w$ is accepted with respect to $\X$ and
$0$ otherwise.  A trivial but useful observation is that if
$A_{\inn}(x, \X_1) = A_{\inn}(x, \X_2)$, then $\sigma(x,
\X_1)=\sigma(x, \X_2)$.

\subsection{Correlation Functions}
\allco

To establish convergence of random point measures to a Gaussian
field, one may employ the method of moments \cite{IS1}, which
depends heavily on the use of correlation functions.

Let $\Po$ denote  $\Po_{\tau}$ or ${\cal Q}$. Recall that
$\alpha(\Po)$ is the subset of $\Po$ which is accepted in the
packing process.  Let $1_B(\alpha(\Po)) = 1$ if $\alpha(\Po) \cap
B$ is non-empty, and let $1_B(\alpha(\Po)) = 0$ otherwise. Given
$w_1,...,w_k \in W$, the $k$-point correlation functions \cite{Ru}
of the point process $\alpha(\Po)$ of accepted
points are defined as
\begin{equation} \label{corr}
r_k(w_1,\ldots, w_k): =\lim_{\epsilon_1 \to 0,\ldots, \epsilon_k
\to 0} \frac {  E_{\Po} 1_{(w_1+ \epsilon_1 \Omega)} \alpha(\Po)
\cdots 1_{(w_k+\epsilon_k \Omega)} \alpha(\Po) }  {
\epsilon_1^{d+1}\cdots \epsilon_k^{d+1}\omega^k},
\end{equation}
where $w_i\neq w_j, 1\leq i\neq j\leq k$, $\Omega$ is the ball in
$W$ of unit radius and  $\omega$ is its volume. To obtain the
definition in the lattice case one replaces the unit balls by the
intervals $\{ {O} \} \times [-1/2, 1/2]$.

To clarify the nature and existence of these correlation functions
we introduce for $w_1,...,w_k \in W$ the functions
$$
\overline{r}_k(w_1,\ldots, w_k):=E_\Po \left[\prod_{i=1}^k
\sigma(w_i,\Po)\right].
$$
The functions $\overline{r}_k$ are the probabilities that all
points  $\{w_1,\ldots, w_k\}$ are packed with respect to an
independent sample from $\Po$. We notice here that this
probability is always positive (unless the balls centered at the
points $x_1,\ldots,x_k$ themselves intersect), regardless of the
heights $t_{x_i}$ of the points.

We want to show that $\overline{r}_k$ are continuous in the
arguments $w_i:=(x_i,t_{x_i}), \ 1 \leq i \leq k$.  We will show
the continuity of  $\overline{r}_k$ in the setting of the
continuum $\Po = \Po_{\tau}$; continuity in the lattice setting
follows from the arguments below and is easier since it only
involves showing continuity in the time coordinate.

Since
$$
\vert \overline{r}_k(w_1,\ldots, w_k) -
\overline{r}_k(w'_1,\ldots, w'_k) \vert
$$
is bounded by
\begin{equation} \label{cont}
\sum_{i=1}^{k-1} \vert \overline{r}_k(w'_1,\ldots,
w'_{i-1},w_i,...,w_k) - \overline{r}_k(w'_1,\ldots, w'_i,w_{i +
1},...,w_k) \vert,
\end{equation}
it will suffice to show that $\overline{r}_k$ is continuous in
each of its $k$ arguments.  We will show
$$
\vert \overline{r}_k(w_1,\ldots, w_k) -
\overline{r}_k(w'_1,\ldots, w_k)\vert \leq C \vert w_1 - w'_1
\vert.
$$
The other summands in (\ref{cont}) are bounded similarly.
 For all $x \in V$, let $B(x)$ denote the ball in
$V$ of unit radius centered at $x$. $B(w), \ w \in W,$ \ is
defined similarly.  Without loss of generality, assume that $\vert
x_1 - x'_1 \vert < 1$, so that $B(x_1)$ and $B(x_2)$ overlap.

Consider the event $\prod_{i=1}^k \sigma(w_i,\Po)\neq\
\sigma(w'_1, \Po) \prod_{i=2}^k \sigma(w_i,\Po)$. If this event
occurs, then
 either the ball $B(x_1)$ is packed and $B(x'_1)$ is not packed or vice versa.   This means
that the oriented graphs on $\Po \cup \{w_1,\ldots,w_k\}$ and $\Po
\cup \{w_1',\ldots,w_k\}$ are different. The only way that this
can happen is if $\Po$ satisfies at least one of the following two
conditions (assume without loss of generality that $t_{x_1} \leq
t_{x'_1}$): \vskip.2cm (a) the cylinder set
$$
[B(x_1) \Delta B(x'_1)] \times [0,t_{x'_1}]
$$
intersects at least one ball $B(w), w \in \Po,$ whereas the
cylinder set
$$
[B(x_1) \cap B(x'_1)] \times [0,t_{x_1}]
$$
does not intersect any ball $B(w), \ w \in \Po$, or \vskip.2cm (b)
the strip
$$
[B(x_1) \cup B(x'_1)] \times [t_{x_1}, t_{x'_1}]
$$
intersects at least one ball $B(w), \ w \in \Po$. \vskip.2cm

Condition (a) implies that

(i) the cylinder centered at $(x_1 + x_1')/2$ with radius $2 -
\vert x_1 - x_1' \vert/2$  and  height $t_{x_1}$ does not contain
any point from $\Po$, whereas

(ii) the cylinder centered at $(x_1 + x_1')/2$  with radius $2 +
\vert x_1 - x_1' \vert/2$  and height $t_{x'_1}$ contains a point
from $\Po$.


Conditions (i) and (ii) imply that there is a set having volume of
order $\vert x_1 - x'_1 \vert$ and  containing at least one point
from $\Po$.  Similarly, condition (b) implies that there is a set
having volume of order $\vert t_{x_1} - t_{x'_1}\vert$ and
containing at least one point from $\Po$. The probability that a
subset of $W$ of small volume contains at least one point from
$\Po$ is of the order of the volume. It follows that the events
(a) and (b) have probability of order $ \vert x_1 - x'_1 \vert +
\vert t_{x_1} - t_{x'_1}\vert$, \ showing that the functions
$\overline{r}_k$ are continuous.






Rewrite the correlation functions $r_k$ as
$$
\lim_{\epsilon_1,\ldots,\epsilon_k\to 0} {E_\Po\left[\prod_{i=1}^k
1_{w_i+\epsilon_i\Omega} (\alpha(\Po)) \ | \ \prod_{i=1}^k
1_{w_i+\epsilon_i\Omega} (\Po) \right]}\times
\frac{E_\Po\left[\prod_{i=1}^k 1_{w_i+\epsilon_i\Omega}
(\Po)\right]}
 { \prod_{i=1}^k \epsilon_i^{d+1} \omega^k}.
$$
We discard the event that any of the small balls
$\{w_i+\epsilon_i\Omega\}_{i = 1}^k$ contains more than one point
of $\Po$ as an event of lower order magnitude. The first factor in
the product is just $\overline{r}_k$ averaged over the possible
positions of these solitary points in the small balls
$w_i+\epsilon_i\Omega$. Continuity of $\overline{r}_k$ implies
that, up to negligible terms, one can replace this factor by the
value of $\overline{r}_k$ at $(w_1,\ldots,w_k)$.
The second factor in the product tends, obviously, to
$$
\prod_{i=1}^k h(w_i),
$$
where $h(w)=h(t)$ is just the density of the underlying Poisson
point process $\Po$ with respect to the volume form on $W$ (in our
situation, $h(t)$ is simply the indicator function of the interval
$[0,\tau]$).

Summarizing, given $w_1,...,w_k \in W$, the correlation functions
are given by the following explicit formula
\begin{equation} \label{corr2}
 r_k(w_1,\ldots, w_k)=\overline{r}_k(w_1,\ldots,
w_k)\prod_{i=1}^k h(w_i).
\end{equation}

If $\alpha(\Po_{\tau})(B)$ represents the number of points in the
set $B \subset \R^d \times \R^+$ which get packed then its moments
are expressed in terms of the one and two point correlation
functions via
\begin{equation} \label{mean}
E[\alpha(\Po_{\tau})(B)] = \int_B r_1(w) dw
\end{equation}
and
$$
E[\alpha^2(\Po_{\tau})(B)] = \int_{B \times B} r_2(w_1,w_2) dw_1
dw_2 + \int_B r_1(w_1)dw_1.
 $$
Also, we have (section 14.4 of Stoyan and Stoyan \cite{SS})
\begin{equation} \label{variance}
\text{Var}[\alpha(\Po_{\tau})(B)] = \int_{B \times B} r_2(w_1,w_2)
- r_1(w_1)r_2(w_2) dw_1 dw_2 +
 \int_B r_1(w_1)dw_1.
\end{equation}


The following is a key definition.

\begin{defn} A family of functions $\{r_j : V^j \to
\R\}_{j=1}^{\infty}$ exponentially clusters if for some positive
constants $A_{k,l}, C_{k,l}$ one has uniformly
\begin{equation*}
|r_{k+l}(w_1,\ldots,w_k,w'_1,\ldots,w'_l)-
r_k(w_1,\ldots,w_k)r_l(w'_1,\ldots,w'_l)| \leq
A_{k,l}\exp(-C_{k,l} d(w,w')),
\end{equation*}
where $d(w,w')$ is the distance between the sets
$\{w_1,\ldots,w_k\}$ and $\{w'_1,\ldots,w'_l\}$, that is the
minimum of the pairwise distances $|w_i-w'_j|$.
\end{defn}

Exponential clustering is also known as weak exponential decrease
of correlations \cite{Ma}, or simply exponential decay of
correlations.


Let $\nu$ be a spatially homogeneous point process on $\R^d$. Then
$\nu$ defines a spatially homogeneous point process $\nu_L$ on the
lattice $\Z^d$ via
$$
\nu_L := \sum_{x \in \Z^d} \nu(Q_x) \delta_x,
$$
where $Q_x$ is the unit cube centered at $x$. Thus $\nu_L$ are
probability measures on $\R^{\Z^d}$.  Exponential clustering of
the correlation functions of $\nu$ implies the exponential
clustering of the correlation functions for $\nu_L$. By following
the methods of Malyshev \cite{Ma}, who restricts attention to
probability measures on $\{-1,1\}^{\Z^d}$, one can show that the
normalized lattice measures
$$
\nu_{L,\lambda}(B):= \frac{\nu_L(\lambda B) - E\nu_L(\lambda B) }
{\sqrt{\lambda^d}}
$$
converge in distribution to a generalized Gaussian field \cite{Ma}
on the Borel sets of $\R^d$ (see Iagolnitzer and Souillard
\cite{IS1}, especially p. 576). The lattice measures
$\nu_{L,\lambda}$ approximate the normalized continuum measures
$$
\nu_{\lambda}(B):= \frac{\nu(\lambda B) - E\nu(\lambda B) }
{\sqrt{\lambda^d}}
$$
in that the mean and variance of $\vert \nu(\lambda B) -
\nu_L(\lambda B) \vert$  are of order $\lambda^{d-1}$ (see
Proposition 4.1 below). Thus the continuum measures
$\nu_{\lambda}$ also converge in distribution to a generalized
Gaussian field.
(An alternative approach involves working directly with the
cumulants and proving that exponential clustering implies
exponential decay of cumulants, automatically yielding the CLT we
are seeking; see \cite{Iv}, \cite{Ba}).

Combining with the expression for the variance (\ref{variance}),
we obtain the following result, the continuum analog of Malyshev's
CLT \cite{Ma} for Gibbsian random fields.


\begin{theo}  \label{FCLT} (Gaussian CLT)
Let $\nu$ be a spatially homogeneous point process such that its
correlation functions exponentially cluster.  Consider the
rescaled centered measures
$$
\nu_{\lambda}(B):= \frac{\nu(\lambda B)-E {\nu} (\lambda B)}
{\sqrt{\lambda^d}}.
$$
Then, as $\lambda \to \infty$,  $\nu_{\lambda}$ converges  in
distribution to a generalized Gaussian random field with
covariance kernel
$$
K(x,y) :=C\delta(x-y),
$$
where
$$
C :=\int_V r_2(0,x)-r_1(0)r_1(x) dx + r_1(0).
$$
\end{theo}

 In other words, for any Borel sets $B_1, \ldots,
B_m$, as $\lambda \to \infty$, \ the vector $\langle
\nu_{\lambda}(B_1),...,\nu_{\lambda}(B_m) \rangle$ tends to a
Gaussian limit with the covariance \ matrix \
$C(\text{vol}(B_i\cap B_j)), \ 1 \leq i, j \leq m$.

\subsection{ Process of packed points}

Recall that $\Po$ denotes either $\P_{\tau}$ or $\Q$. To show
Theorems 1.2 and 1.3, we will show that the correlation functions
$r_k, \ k \geq 1,$ of the point process $\alpha(\Po)$ cluster
exponentially. This will imply that the correlation functions
$r_k^{\pi}$ of the point process $\pi(\Po)$ cluster exponentially,
which, by  Theorem \ref{FCLT}, gives the desired result.

To show exponential clustering of the $r_k, \ k \geq 1,$ we first
establish that $\alpha(\Po)$ is localized near $V$. This
localization is of course obvious in the off-lattice case, as the
$t$-support of the process is the bounded interval $[0,\tau]$. We
cannot remove this cut-off as the current way of proof needs a
rather rapid decay of the correlation functions with $t$ which is
lacking in the off-lattice case. Indeed, we have some basic
estimates on bounds for the correlation functions.  For example,
the decay of $r_1$ is polynomial: in dimension $d=1$ it follows
directly from the R\'enyi formula that $r_1(x,t)\sim t^{-2}$ and
in higher dimensions $r_1$ can be also shown to be at least of
order $t^{-\frac{1}{d}-1}$. The next lemma shows that the
correlation functions decay exponentially in $t$ in the finite
input off-lattice case as well as in the infinite input lattice
setting.

\begin{prop}\label{A} The correlation functions $r_k, \ k \geq 1,$ of the
point process $\alpha(\P)$  decay exponentially with $t$'s, i.e.
for any $k$ there exists positive constants $A_k$ and $C_k$ such
that for $w_1,...,w_k \in W$ we have
$$
r_k(w_1,\ldots,w_k)\leq A_k\exp(-C_k \max_i(t_i)).
$$
\end{prop}

{\em Proof.} Let $w_i = (x_i, t_i), \ 1 \leq i \leq k$.
 Consider first the correlation functions for
$\alpha(\Po_{\tau})$.  Proposition \ref{A} is clearly satisfied
 as in this case the left hand side of the
inequality vanishes for $ \text{max}_{i \leq k} t_i \geq \tau$.

Now consider the correlation functions for $\alpha({\cal Q})$,
 that is the correlation functions for the
point process of accepted points on the lattice. Notice that,
obviously,
$$
\overline{r}_k(w_1,\ldots,w_k)\leq \overline{r}_1(w_i)
$$
for any $1\leq i\leq k$. Further, $\overline{r}_1(w_1)$ is just
the probability that the point $w_1=(x,t)$ is packed with respect
to ${\cal Q}$. However, this implies that all points $w_s=(x,st),
0<s<1,$ can be packed with respect to ${\cal Q}$, and therefore,
that none of the points $w_s, \ s<1,$ are present in the sample
${\cal Q}$. This has probability $\exp(-a t)$, where $a >0$ is the
intensity measure of the interval $\{(x,s), \ 0\leq s\leq 1\}$.
Thus $r_1=\overline{r}_1$ decays exponentially with $t$ and
Proposition \ref{A} is proved in the lattice setting. \qed

\vskip.3cm

\begin{prop} \label{CC} Let $E_{k,l}:= E_{k,l}(w_1,...,w_k,w_1',...,w_l';\Po)$
be the event that the
backward cones of any pair of points $w_i$ and $w_j'$
with respect to $\Po$ do not intersect. Then
\begin{equation}\label{3.2}
\vert r_{k + l}(w_1,...,w_k,w_1',...,w_l') - r_{k} (w_1,...,w_k)
r_{l} (w_1',...,w_l') \vert \leq C_{k,l} (P[E^c_{k,l}])^{1/2}.
\end{equation}
\end{prop}

{\em Proof.} Define for any event $E$ measurable with respect to
the sigma algebra generated by $\Po$
$$
\overline{r}_k^E(w_1,\ldots, w_k): =
E_\Po \left[\prod_{i=1}^k
\sigma(w_i,\Po) 1_E\right].
$$
It suffices to show that (\ref{3.2}) holds with $r_k$ replaced by
$\overline{r}_k$.
Now $|\overline{r}_k - \overline{r}_k^E|$ equals
\begin{equation} \label{3.3}
E_\Po \left[\prod_{i=1}^k \sigma(w_i,\Po) 1_{E^c}\right] \leq A_k
(P[E^c])^{1/2},
\end{equation}
by Cauchy-Schwarz and the boundedness of $\overline{r}_k$.

Observe that by the definition of $E_{k,l}$, the random variables
$$
\prod_{i=1}^k \sigma(w_i,\Po)1_{E_{k,l}} \ \ \text{and} \ \
\prod_{j=1}^l \sigma(w_j',\Po)1_{E_{k,l}}
$$
are independent.

Hence the difference $\overline{r}^{E_{k,l}}_{k +
l}(w_1,...,w_k,w_1',...,w_l') - \overline{r}^{E_{k,l}}_{k}
(w_1,...,w_k) \overline{r}^{E_{k,l}}_{l} (w_1',...,w_l')$ vanishes
and using the estimate (\ref{3.3}) we obtain (\ref{3.2}). \qed

The clustering of the correlation
functions is captured in

\begin{prop}\label{B} The correlation functions $r_k, \ k \geq 1,$
of the point process $\alpha(\P)$ cluster exponentially.
\end{prop}

{\em Proof.} We fix $\tau$ and use Lemma \ref{cluster} describing
the localization of  causal cones. Let $\P = \P_{\tau}$;  the
proof for $\P = \Q$ is exactly the same. The constants $\gamma,
\rho,$ and $\beta$ are as in Lemma \ref{cluster}. Let
$w_i:=(x_i,t_i), \ 1 \leq i \leq k$. Let $w_j' := (x_j',t_j'), \ 1
\leq j \leq l$.

Consider two sets $\{w_1,\ldots,w_k\}$ and $\{w'_1,\ldots,w'_l\}$
at distance $d$. That is $d := \min |w_i - w_j' |$.
 We distinguish two cases.

(a) All times $t_1,\ldots,t_k, t'_1,\ldots,t'_l$ are less than
$\frac{d}{4\beta}$. In this case with probability at least $1- (k
\rho)\exp(-\gamma d/4)$, the causal cone for each point $w_i, \ 1
\leq i \leq k,$ belongs to the set $C_{d/4}(w_i)$ which is a
subset of the cylinder of radius $d/2$ centered at $x_i$. The same
is valid for $w'_j, \ 1 \leq j \leq l$. Let $\E$ denote the event
for which this is true. As corresponding cylinders for points from
different tuples do not intersect, we conclude that $\E$ implies
$E_{k,l}$, where $E_{k,l}$ is defined as in Proposition \ref{CC}.
Now apply Proposition \ref{CC}.

(b) If at least one of the times (say, $t_1$) is larger than
$\frac{d}{4\beta}$, then by Proposition \ref{A}, both
$r_{k+l}(w_1,\ldots w_k, w_1',...,w_l')$ and $r_k(w_1,\ldots,w_k)$
decay exponentially in $t_1$ and therefore decay exponentially in
$d$ as well. This finishes the proof. \qed

\vskip.5cm

We want to show the Gaussian structure  for the process $\pi(\Po)$
by applying Theorem \ref{FCLT}. $\pi(\Po)$ is a translationally
invariant discrete (meaning that the distance between any two
points is uniformly bounded from below) point process.

 To this end, we have
to check that the correlation functions cluster exponentially.
This is easy, given Propositions \ref{A} and \ref{B}. Indeed,
given $x_1,...,x_k \in V$,  the correlation functions for
$\pi(\Po)$ are given as
$$
r^{\pi}_k(x_1,\ldots,x_k):= \int_{t_1,\ldots,t_k\geq 0}
r_k((x_1,t_1),\ldots, (x_k,t_k)) dt_1\ldots dt_k.
$$

 To prove the
clustering inequality for a $k$-tuple $\langle x_1,...,x_k\rangle$
and an $l$-tuple $\langle x_1',...,x_l'\rangle$ at a distance $d$
we split the integration domain $(\R^+)^k \times (\R^+)^l$ into
two subdomains, namely $[0,d]^{k + l}$ and its complement. Since
all $t_i, \ 1 \leq i \leq k,$ and $t'_j, \ 1 \leq j \leq k,$ are
less than $d$ in the first domain, we apply Proposition \ref{B},
using the fact that the $W$-distance between $w_i$ and $w'_j$ is
at least $d$ and the polynomial bound $d^{k + l}$ on the volume of
the integration domain.  This gives a bound which is exponentially
decaying with $d$.

In the second subdomain we apply Proposition \ref{A}, resulting in
the estimate from above of this second part as
$$
\int_d^\infty A_{V,k,l}\exp(-C_{V,k,l}s) d
\text{vol}(\{\max(t_1,\ldots,t_k,t'_1,\ldots,t'_l)\leq s\}),
$$
where $A_{V,k,l}:=A_k+A_l+A_{k+l}$ and $C_{V,k,l}:=\min(C_k, C_l,
C_{k+l})$. This is obviously exponentially decaying with $d$ and
thus the correlation functions $r^{\pi}_k$ cluster exponentially.

Now applying Theorem \ref{FCLT} we obtain the Gaussian CLT for the
rescaled measures $\nu_{\tau,\lambda }$. This completes the proof
of Theorem \ref{main1}. \qed

\vskip1cm

\section{Proof of Theorem \ref{main2}}
\allco

 We now deduce Theorem \ref{main2} from Theorem
 \ref{main1}.  We have already noted that because of edge effects,
 the rescaled point measures
 $\nu_{\tau,\lambda}$ and $\nu_{\tau,A, \lambda}$ are in general
 not equal.
However, one can estimate the difference between the point
measures $\nu_{\tau,\lambda}$ and $\nu_{\tau, A, \lambda}$.  We
will do this by using  the exponentially clustering of the
two-point correlation function to upper bound the  variance of the
difference.

More specifically, for an open set $A \subset \R^d,$ let
$\pi(\Po)_A := \pi(\Po) \cap A$ and define the point measure
$\pi^+_{\lambda A}:= \sum_{x \in \pi(\Po_{\lambda A})-
\pi(\Po)_{\lambda A}} \delta_x $. $\pi^+_{\lambda A}$ is the
difference between the packing process on $\lambda A$ and the
infinite packing process restricted to $\lambda A$. Similarly
define
 $\pi^-_{\lambda A} := \sum_{x \in \pi(\Po)_{\lambda A}- \pi( \Po_{\lambda A}) } \delta_x
 $.

We want to estimate the variance of the number of points defining
the supports of the point measures $\pi^+_{\lambda A}$ and
$\pi^-_{\lambda A}$. If we could show that this variance is of
order  $o(\lambda^d),\ \lambda \to\infty$, then it would follow
that the random measures $\pi(\Po_{\lambda A})$, after centering
and rescaling by $\sqrt{\lambda^d},$ have asymptotically the same
distribution as the centered and rescaled random measures
$\pi(\Po)_{\lambda A}$, that is, have a generalized Gaussian
distribution.

To prove these variance bounds, we again resort to representation
of the moments as integrals of polynomials of the correlation
functions. In our case we will need just the expressions
(\ref{mean}, \ref{variance}) for the first two moments.

For all $k$ we let $r_k^\pm$ denote the correlation functions for
$\pi^\pm_{\lambda A}$.  Observe the following:

(a) The correlation function $r'_1(w):= P[ w\in \pi(\Po_{\lambda
A} )| w\in \Po_{\lambda A}]$ decays exponentially with time $t$
for the same reasons that $r_1$ does.

(b) The correlation functions $r_k^\pm$ cluster exponentially. The
proof goes along the lines of the proof of Proposition \ref{B},
and in fact involves conditioning on the same event.

In other words, Propositions \ref{A} and \ref{B}, together with
the constants there, are valid for the point processes
$\pi_{\lambda A}^\pm$.

 The crucial property is the next one, which
ensures that the point processes $\pi^\pm_{\lambda A}$ have
support which localizes near the boundary of $\lambda A$.

(c)
The causal cone $A_{\out,\inn}(w, \Po_{\lambda A})$ of any point
$w$ with respect to the process $\Po_{\lambda A}$ is a subset of
$A_{\out, \inn}(w, \Po_{\lambda A})$ with respect to $\Po$ and
coincides with it if $A_{\out, \inn}(w, \Po)$ does not intersect
the boundary $\partial(\lambda A)$.

Therefore, the correlation functions $r_1^\pm$ of the processes
$\pi^\pm_{\lambda A} $ decay exponentially with the distance of
$V$-projection to $\partial(\lambda A)$:
$$
r_1^\pm(w)\leq A_\partial \exp\{-C_\partial d(x, \partial(\lambda
A))\}, \ \ w := (x,t)
$$
for some positive constants $A_\partial,C_\partial$.

Together these properties imply

\begin{prop} \label{C} Assume that the set $A$ is piecewise smooth.
Then the variance and the mean of the number of points in either
of the processes $\pi^\pm_{\lambda A}$ is of order
$O(\lambda^{d-1})$.
\end{prop}

\vskip.5cm

{\em Proof. } Start with the mean, which by (\ref{mean}), is given
by
$$
\int_W r^\pm_1(w) dw.
$$
By properties (a) and (c) above, the correlation function
$r_1^\pm$ decays exponentially (say, as $D\exp(-Cd)$) with the
distance $d$ from $\partial(\lambda A)\subset V\times\{0\}\subset
W$. Hence the integral in question is bounded above by
$$
D \int_0^\infty \exp(-C s) dB_\lambda(s),
$$
where $B_\lambda(s):=\text{vol}(\{x \in \R^d: \ d(x,
\partial(\lambda A))\leq s\})$ and where $\{x \in \R^d:\ d(x,
\partial(\lambda A))\leq s\}$ is the $s$-neighborhood around the
set $\lambda A$. Since the leading term for the volume of the
$s$-neighborhood around $\lambda A$ has the form $s \lambda^{d-1}
C_A$, were $C_A$ is a constant depending on the set $A$, the above
integral has leading term in $\lambda$ which is polynomial of
degree $(d-1)$.

To estimate the variance, we just add the clustering condition and
use the variance formula (\ref{variance}) to obtain
$$
\text{var}(\pi^+_{\lambda A})= \int_{W\times W} [r^+_2(w_1,
w_2)-r^+_1(w_1)r^+_1(w_2)] dw_1 dw_2+ \int_W r^+_1(w) dw,
$$
where $r^+_1$ and $r^+_2$ are the correlation functions for the
process $\pi^+_{\lambda A}$ (an analogous formula is also valid
for $\text{var}( \pi^-_{\lambda A})$).

The second summand is just the mean, and the integrand in the
first integral decays in $W\times W$ exponentially with the
distance from $\partial(\lambda A)\subset W\times W$. Once again
applying the integral estimates as above, we obtain
$\text{var}(\pi^+_{\lambda A}) = O(\lambda^{d-1}).$ \qed

\section{Gaussian Fields and Total Edge Length Functionals}
\allco

The above methods extend existing central limit theorems involving
functionals of Euclidean point sets, including those in
computational geometry, Euclidean combinatorial optimization, and
Boolean models. Existing central limit theorems (\cite{PY1} and
references therein) show asymptotic normality of functionals such
as total edge length, total number of components, and total number
of vertices of a specified degree.

These functionals are canonically associated with point measures
on $\R^d$.  Given a graph $G$ on a locally finite point set $\X$,
we associate to the total edge length functional  the point
measure  defined by giving each vertex $x \in \X$ a weight equal
to one half the length of the edges in $G$ incident to $x$.


For example, given $\X$ a point set, let $\nu_{\X}:= \nu_{NN(\X)}$
be the point measure associated with the total edge length of the
nearest neighbors graph on $\X$. Thus, if $G(\X)$ is the nearest
neighbors graph on $\X$ and if $\E(x;G(\X))$ denotes the edges
incident to $x,$ then $\nu_{NN(\X)}$ is defined as
$$
\nu_{NN(\X)}(B):= \sum_{x \in B \cap \X} \frac{1}{2} \sum_{e \in
\E(x; G(\X)) } \vert e \vert,
$$
where $\vert e \vert$ denotes the length of the edge $e \in G$.
Point measures associated with the total edge length of the
Voronoi tessellation, minimal spanning tree, and sphere of
influence functionals, are defined analogously.

Let $\nu_{NN}:= \nu_{NN,\P}$ (respectively, $\nu_{NN,A}:=
\nu_{NN,A,\P}$) denote the random measures associated with the
nearest neighbors graph on the Poisson point process $\P$
(respectively, $\P \cap A$). Analogously to the rescaled packing
measures (\ref{measure}) and (\ref{measure1}),  define the
rescaled (infinite volume) measures:
\begin{equation} \label{measure3}
\nu_{NN, \lambda}(B):= \frac{\nu_{NN} (\lambda B)- E \nu_{NN}
(\lambda B)} {\sqrt{\lambda^d}}
\end{equation}
and the rescaled (finite volume)  measures
\begin{equation}\label{measure4}
\nu_{NN, A, \lambda}(B):= \frac{\nu_{NN, \lambda A} (\lambda B)- E
\nu_{NN, \lambda A} (\lambda B)} {\sqrt{\lambda^d}}.
\end{equation}

Consider now the correlation functions for the point measure
$\nu_{NN(\X)}$. Define the correlation functions of the point
measures $\nu_{NN}$ by:
$$
r_k(w_1,\ldots, w_k): =\lim_{\epsilon_1 \to 0,\ldots, \epsilon_k
\to 0} \frac {  E_{\Po} \nu_{NN}(w_1+ \epsilon_1 \Omega) \cdots
\nu_{NN}(w_k+\epsilon_k \Omega)) }  { \epsilon_1^{d+1}\cdots
\epsilon_k^{d+1}\omega^k},
$$

In contrast to the packing measures, the point measures $\nu_{NN}$
do not evolve with time $\tau$.  Instead of considering their
localization properties via causal cones in $\R^d \times
[0,\tau]$, we only need to consider their localization properties
in $\R^d$.  This is accomplished by considering the notion of
stabilizing measures \cite{PY2}.
Given a point measure $\nu$, write $\nu_{\X}$ when $\nu$ has
support $\X$ on $\R^d$. Say that $\nu$ is {\em stabilizing
exponentially fast} if there is a ball $B_R: = B_R(0)$, centered
at the origin, with an exponentially decaying radius $R$, that is,
$P[R > t] \leq \exp(-Ct)$, such that the difference measures
defined by
$$
\nu_{\Po \cap B_R \cup {\cal A} \cup {\bf 0} }(\ \cdot \ ) -
\nu_{\Po \cap B_R \cup {\cal A} }(\ \cdot \ )
$$
are invariant for all finite ${\cal A} \subset \R^d - B_R(0).$ In
other words, changes in the environment outside $B_R(0)$ do not
change the values of the point measure $\nu$  inside $B_R(0)$.

By straightforward modifications to the proof of Proposition 3.3,
it follows that if point measures stabilize exponentially fast,
then their correlation functions will cluster exponentially, and
therefore, the rescaled point measures converge to a Gaussian
field.

Now we may show \cite{PY2} that the point measure $\nu_{NN}$
stabilizes exponentially fast.  We thus have the following
convergence result.  An analogous convergence result holds for the
rescaled finite volume measures (\ref{measure4}).

\begin{theo} (Infinite volume nearest neighbor measures converge to a Gaussian field)
\label{main4} The measures $\nu_{NN, \lambda}$ converge in
distribution as $\lambda \to \infty$ to a generalized Gaussian
random field with covariance kernel $K$ concentrated on the
diagonal, that is
$$
K (x,y) :=C \delta(x-y),
$$
where $C$ is a constant.
\end{theo}

Moreover, if  we consider Voronoi tessellations, sphere of
influence graphs, or minimal spanning tree graphs on Poisson point
sets, then the canonically associated point measures associated
with the total edge length functional are exponentially
stabilizing \cite{PY2} and they thus converge to a generalized
Gaussian random field.

\end{document}